\documentclass[11pt]{article}
\usepackage[T1]{fontenc}
\usepackage{lmodern,microtype}
\usepackage[margin=1in]{geometry}
\usepackage{tikz-cd}
\usepackage{amsmath,amssymb,amsthm,mathtools}
\usepackage{thmtools}
\usepackage{enumitem,booktabs,array,xcolor,titlesec,fancyhdr}
\definecolor{linkblue}{HTML}{24527A}
\usepackage[colorlinks=true,linkcolor=linkblue,citecolor=linkblue,urlcolor=linkblue]{hyperref}
\hypersetup{pdftitle={On the Topology of the Set of Nonnegative Forms That Are Not Sums of Squares},pdfauthor={}}
\usepackage[nameinlink,noabbrev]{cleveref}
\newcommand{\partlabel}[1]{\textup{(}\ensuremath{#1}\textup{)}}
\setlist[enumerate]{label=\partlabel{\roman*},leftmargin=*,itemsep=2pt,topsep=4pt,parsep=0pt}
\titleformat{\section}{\large\bfseries}{\thesection}{0.75em}{}
\titleformat{\subsection}{\normalsize\bfseries}{\thesubsection}{0.75em}{}
\titlespacing*{\section}{0pt}{2.5ex plus 0.6ex minus 0.2ex}{1.1ex}
\titlespacing*{\subsection}{0pt}{2.0ex plus 0.4ex minus 0.2ex}{0.8ex}
\allowdisplaybreaks[1]
\numberwithin{equation}{section}

\newtheorem{theorem}{Theorem}[section]
\newtheorem{proposition}[theorem]{Proposition}
\newtheorem{lemma}[theorem]{Lemma}
\newtheorem{corollary}[theorem]{Corollary}
\theoremstyle{definition}

\theoremstyle{remark}

\newcommand{\R}{\mathbb R}
\newcommand{\Q}{\mathbb Q}
\newcommand{\Z}{\mathbb Z}

\newcommand{\RP}{\mathbb {RP}}
\newcommand{\Pos}{\operatorname{Pos}}
\newcommand{\Sos}{\operatorname{Sos}}
\newcommand{\Sym}{\operatorname{Sym}}
\newcommand{\tr}{\operatorname{tr}}

\newcommand{\Span}{\operatorname{span}}
\newcommand{\cone}{\operatorname{cone}}
\newcommand{\interior}{\operatorname{int}}
\newcommand{\Conf}{\operatorname{Conf}}
\newcommand{\UConf}{\operatorname{UConf}}
\newcommand{\SO}{\operatorname{SO}}
\newcommand{\eps}{\varepsilon}
\newcommand{\tH}{\widetilde H}
\newcommand{\tchi}{\widetilde\chi}

\title{On the Topology of the Set of Nonnegative Forms That Are Not Sums of Squares}
\author{Maximilian Illmer \quad and \quad Tim Netzer}
\date{\today}

\makeatletter
\renewcommand{\maketitle}{%
  \begin{center}
    {\LARGE\bfseries \@title\par}
    \vspace{0.8em}
    {\large \@author\par}
    \vspace{0.4em}
    {\normalsize Department of Mathematics, University of Innsbruck\par}
    \vspace{0.6em}
    {\small \@date\par}
  \end{center}
  \vspace{1.2em}
  \thispagestyle{plain}%
}
\makeatother

\renewenvironment{abstract}{%
  \begin{center}
  \begin{minipage}{0.88\textwidth}
  \small\noindent\textbf{Abstract.}\enspace
}{%
  \end{minipage}
  \end{center}
  \vspace{0.8em}
}

\begin{document}
\maketitle
\begin{abstract}
We study the topology of the normalized space of real nonnegative
forms of degree $2d$ in $n$ variables that are not sums of squares.
Whenever it is nonempty, we prove that it is path-connected but not acyclic, and in particular not contractible. We give explicit bounds on the degrees in which its reduced homology can be
nonzero. In the two minimal
nonempty cases, we show an embedding of an exterior algebra into the cohomology ring, and prove containment of the binary octahedral and binary icosahedral
groups, respectively, in its fundamental group.
\end{abstract}

\section{Introduction}
A basic problem in real algebra is to relate global nonnegativity of polynomials, a condition on all real points, to algebraic representations such as sums of squares. The two notions agree in some cases but not in general. Hilbert's 1888 theorem gives an exact classification: every nonnegative homogeneous form of positive even degree is a sum of squares only in one of the following three cases of families: quadratic forms in any number of variables, forms in two variables, and ternary quartics \cite{Hilbert}. A decade later, Hilbert's 17th problem asked whether every polynomial that is nonnegative on $\R^n$ can be expressed as a sum of squares of rational functions \cite{HilbertProblems}. Artin answered this affirmatively in 1927 \cite{Artin}. 
 
Over the following decades, work on ordered and real-closed fields and quantifier elimination, together with the real Nullstellensatz and Positivstellens\"atze, developed a systematic framework for describing real solutions and certifying nonnegativity, see \cite{BCR,Marshall,PowersBook,PrestelDelzell,ScheidererCourse}. In modern work, sums of squares are also used as tractable certificates in polynomial optimization: Gram matrix representations turn their search into semidefinite optimization; see 
for example \cite{LaurentSOS}.

Let $n\geqslant2$, $d\geqslant1$, and let $H_{n,d}$ denote the real vector space of homogeneous polynomials of degree $d$ in $n$ variables. Throughout,
\[
 m=\dim H_{n,2d}=\binom{n+2d-1}{2d},\qquad
 k=\dim H_{n,d}=\binom{n+d-1}{d},\qquad
 c=m-nk+\binom n2-2.
\]
Write $\Pos_{n,2d}$ and $\Sos_{n,2d}$ for the convex cones of globally nonnegative forms and sums of squares of degree-$d$ forms. To obtain compact convex sets we normalize by
\[
 L(f)=\int_{S^{n-1}} f\,d\sigma,
 \qquad P=\Pos_{n,2d}\cap \{L=1\},\qquad
 S=\Sos_{n,2d}\cap\{L=1\},
\]
where $\sigma$ is the rotation-invariant probability measure. Our main object of study is the set
\[
 E=P\setminus S,
\]
where our aim is to understand its topology.
We write $E_{n,2d}$ for the set $E=P\setminus S$ when we wish to emphasize its dependence on $n$ and $2d$.
In particular, $E_{3,6}$ and $E_{4,4}$ denote the sets of normalized
nonnegative ternary sextics and quaternary quartics, respectively, that are
not sums of squares.
Hilbert's theorem shows that $E$ is empty exactly when $d=1$, $n=2$, or $(n,2d)=(3,4)$. We exclude these cases unless stated otherwise.

Although Hilbert's theorem showed that nonnegative forms are not always sums of squares, the first explicit example appeared only in Motzkin's 1967 work, nearly eighty years later; Robinson subsequently gave further examples \cite{Motzkin,Robinson}. Against this background, Blekherman's volume estimate is especially striking: it quantifies the difference between the sets $S$ and $P$ \cite{BlekhermanVolume}. His result implies that for each fixed $d\geqslant2$ there is a constant  $C_d>0$, depending only on $d$, such that 
\[
 \left(\frac{\operatorname{vol}(S)}{\operatorname{vol}(P)}\right)^{1/(m-1)}
 \leqslant C_d n^{-(d-1)/2}.
\]
Here, $\operatorname{vol}$ denotes $(m-1)$-dimensional Lebesgue measure on the hyperplane $\{L=1\}$. Thus, after taking the appropriate root, the sum-of-squares section becomes small relative to the nonnegative section, as the number of variables grows. 

To the best of our knowledge, the \emph{topology} of $E$ has not previously been studied. The volume estimates above quantify the relative size of $S$ inside $P$, but leave open the shape of its complement $E$; neither those estimates nor the convexity of $P$ and $S$ determines that shape. We show that whenever $E$ is nonempty, it is path-connected but not contractible, has nonzero reduced rational homology, and is homotopy equivalent to a finite CW complex with an explicit dimension bound. We also locate all possible nonzero reduced homology groups in an explicit degree range. In broad parameter ranges its Euler characteristic vanishes, forcing odd-degree homology. The two minimal nonempty cases exhibit still richer structure: exterior subalgebras in rational cohomology, direct summands in $\pi_3$, and binary octahedral or icosahedral subgroups in $\pi_1$. Together these results give a systematic topological picture of the gap between nonnegativity and sums of squares.

We state the results precisely below. The first theorem describes the topology in every nonempty case, while the second gives sharper information in the two minimal cases. All homology groups without an indicated coefficient ring have coefficients in $\Z$.

\begin{theorem}\label{thm:main-general}
Whenever $E$ is nonempty, the following assertions hold.
\begin{enumerate}
\item The space $E$ is path-connected and  homotopy equivalent to a finite CW complex of dimension at most $m-k-1$.
\item All nonzero reduced integral homology of $E$ lies in the range
\begin{equation}\label{eq:range}
 \max\{1,c+1\}\leqslant i\leqslant m-k-1.
\end{equation}

\item The Euler characteristic $\chi(E)$ is even. In particular, $E$ has nonzero reduced rational homology and is not contractible. Moreover,
\[
 \chi(E)=0\quad\text{if}\quad
 \begin{cases}
 n=3,\ d\geqslant3,\quad\text{or}\quad n=4,\ d\geqslant2,\\
 \text{or}\quad 5\leqslant n\leqslant8,\ d=2.
 \end{cases}
\]
In these cases some odd-degree rational homology group in the range from \eqref{eq:range} is nonzero.
\end{enumerate}
\end{theorem}

In the two minimal nonempty cases, we obtain additional information.
Here and below, $\Lambda_{\Q}(u_3,\ldots)$ denotes the exterior algebra over
$\Q$, with subscripts indicating the degrees of the generators. We write $b_i$ for the Betti number of dimension $i$.

\begin{theorem}\label{thm:main-small}
$(i)$ For ternary sextics, there is an embedding of graded algebras
\[
 \Lambda_{\Q}(u_3)\hookrightarrow H^*(E_{3,6};\Q),
\]
which implies $b_3(E_{3,6})\geqslant1$. Furthermore, $\pi_3(E_{3,6})$ contains a direct summand isomorphic to $\mathbb Z$, $\pi_1(E_{3,6})$ contains a subgroup isomorphic to  the binary octahedral
group $2O$ of order $48$, and the center of $\pi_1(E_{3,6})$ contains a copy of $\Z/2$.

$(ii)$ For quaternary quartics, there is an embedding of graded algebras
\[
 \Lambda_{\Q}(u_3,v_3)\hookrightarrow H^*(E_{4,4};\Q),
\]
which implies $b_3(E_{4,4})\geqslant2$ and $b_6(E_{4,4})\geqslant1$.
Furthermore, $\pi_3(E_{4,4})$ contains a direct summand isomorphic to $\Z^2$, $\pi_1(E_{4,4})$ contains the binary icosahedral group $2I$ of order $120$, and the center of $\pi_1(E_{4,4})$  contains a copy of $\Z/2$.
\end{theorem}

In Sections 2--5 we establish the results that hold in all cases, and Section 6 treats the two minimal cases.

\section{Convex geometry and basic homotopies}
Consider the auxiliary sets
\[
 \interior{E}=\interior(P)\setminus S,
 \qquad F\coloneqq S\cap\partial P,
\]
where interiors and boundaries are taken in the affine hyperplane $\{L=1\}$.
For $h\in P$, let $Z(h)\subseteq\RP^{n-1}$ be its real projective zero set and let $I_r(Z(h))$ be the space of real degree-$r$ forms vanishing there.
We also use the set
\[
 D\coloneqq\{h\in P\mid I_d(Z(h))=0\}.
\]
Every member of $D$ lies in $E$. Indeed, if $h=\sum_i q_i^2$, then each $q_i$ vanishes on $Z(h)$ and hence belongs to $I_d(Z(h))=0$, contradicting $L(h)=1$.

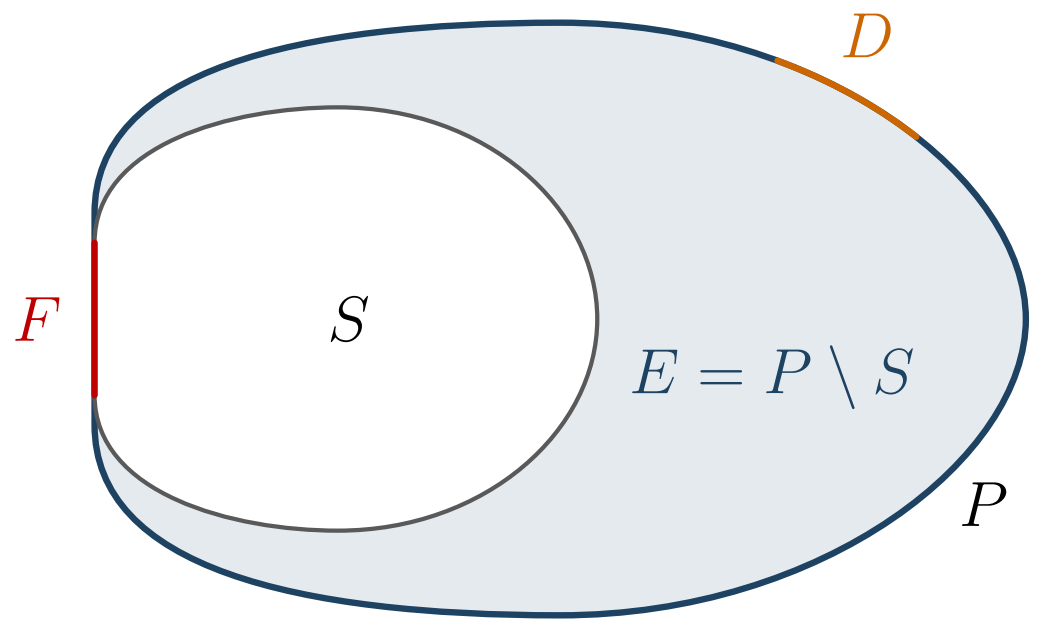
\begin{figure}[htbp]
\centering
\begin{tikzpicture}[scale=0.70, line cap=round, line join=round]
  \path[fill=linkblue!12,draw=linkblue!80!black,line width=1.2pt]
    (-4.4,-1.0) -- (-4.4,1.0)
    .. controls (-4.4,2.35) and (-2.5,2.8) .. (0,2.8)
    .. controls (2.5,2.8) and (4.4,1.3) .. (4.4,0)
    .. controls (4.4,-1.3) and (2.5,-2.8) .. (0,-2.8)
    .. controls (-2.5,-2.8) and (-4.4,-2.35) .. (-4.4,-1.0) -- cycle;
  \path[fill=white,draw=black!65,line width=0.8pt]
    (-4.4,-0.72) -- (-4.4,0.72)
    .. controls (-4.4,1.65) and (-3.1,2.0) .. (-2.1,2.0)
    .. controls (-0.7,2.0) and (0.35,1.05) .. (0.35,0)
    .. controls (0.35,-1.05) and (-0.7,-2.0) .. (-2.1,-2.0)
    .. controls (-3.1,-2.0) and (-4.4,-1.65) .. (-4.4,-0.72) -- cycle;
  \draw[very thick,red!75!black] (-4.4,-0.72) -- (-4.4,0.72);
    \draw[very thick,orange!80!black] (2.0529,2.4409)
    .. controls (2.5587,2.2542) and (3.0025,2.0055) .. (3.3642,1.7216);
  \node[font=\large] at (-2.0,0) {$S$};
  \node[font=\large,text=linkblue!80!black] at (2,-0.55) {$E=P\setminus S$};
  \node[font=\large] at (4,-1.75) {$P$};
  \node[font=\large,text=orange!80!black] at (2.9,2.68) {$D$};
  \node[font=\large,text=red!75!black,anchor=east] at (-4.52,0) {$F$};
\end{tikzpicture}
\caption{\small Schematic of the involved sets: $P$ is the whole body, $S\subseteq P$ is the sum-of-squares region, $E=P\setminus S$ is shaded, and $F=S\cap\partial P$ is shown as a shared boundary segment. The orange arc indicates $D\subseteq\partial P\setminus F\subseteq E$. This is also close to the picture we had in mind when starting this project, making the nontrivial topology of $E$ that we discovered all the more suprising.}
\label{fig:set-relations}
\end{figure}

Fix an orthonormal basis
$
 v=(v_1,\ldots,v_k)^T
$
of $H_{n,d}$ for the inner product $\langle q,r\rangle=L(qr)$, and put
\[
 \Gamma(Q)=v^TQv\quad\mbox{ for}  Q\in\Sym_k(\R).
\]
Then $L(\Gamma(Q))=\tr Q$. This Gram map is onto, since every degree-$2d$ monomial is a product of degree-$d$ monomials. Its image of the positive definite matrices is open and lies in the sos cone.

Quantifier elimination shows that all sets under consideration are semialgebraic. The trace identity also shows that the sos cone is closed: Gram matrices of a convergent sequence of sos forms have bounded trace, so a subsequence converges. Since $L$ is strictly positive on every nonzero nonnegative form, $P$ is compact; otherwise normalized unbounded forms would yield a nonzero nonnegative form with integral zero. Thus $P$ and $S$ are compact convex bodies in $\{L=1\}\cong\R^{m-1}$. The interior of $P$ consists exactly of strictly positive forms. It follows that $\interior{E}$ is nonempty in every nonempty case.

We first establish some basic homotopies that allow us to apply Alexander duality. This only uses convexity of the involved sets.

\begin{proposition}\label{prop:collar}
Choose $c_0\in\interior(S)$, and let $r_S(u)$ and $r_P(u)$ be the radial functions of $S$ and $P$ about $c_0$, for $u\in S^{m-2}$. Set
\[
 \Omega=\{u\mid r_S(u)<r_P(u)\}.
\]
Then there are homeomorphisms
\[
 E\cong\Omega\times(0,1],\qquad \interior{E}\cong\Omega\times(0,1),\qquad
 \partial P\setminus F\cong\Omega.
\]
In particular $E$ is an orientable topological $(m-1)$-manifold with boundary, $\interior{E}\hookrightarrow E$ is a homotopy equivalence, and radial projection is a strong deformation retraction from $E$ onto $\partial P\setminus F$. Furthermore,
\begin{equation}\label{eq:alexander}
 \tH_i(E)\cong\tH^{m-i-3}(F).
\end{equation}
\end{proposition}
\begin{proof}
The homeomorphism $\Phi\colon\Omega\times(0,1]\to E$ is given by
\[
 \Phi(u,t)=c_0+\bigl((1-t)r_S(u)+tr_P(u)\bigr)u.
\]
Its inverse sends $x\in E$ to
\[
 u=\frac{x-c_0}{\|x-c_0\|},
 \qquad
 t=\frac{\|x-c_0\|-r_S(u)}{r_P(u)-r_S(u)}.
\]
Since the radial functions are continuous and
$r_P(u)-r_S(u)>0$ on $\Omega$, continuity of both $\Phi$ and its inverse is
readily checked.  Restricting the same formulas gives the other asserted
homeomorphisms and the retraction. For the co-/homology statement we use that $F$ is compact and locally contractible by semialgebraicity, so that Alexander duality on $\partial P\cong S^{m-2}$ gives \eqref{eq:alexander}.
\end{proof}

For a nonnegative form $h$, define its contact cone by
\[
 C_h=\cone\{v(a)v(a)^T\mid [a]\in Z(h),\ \|a\|=1\},
 \qquad \rho(h)=\dim C_h.
\]
The nonzero generators have traces bounded away from zero, so $C_h$ is
closed. For normalized $h$ one has
\[
 h\in D\quad\Leftrightarrow\quad C_h\text{ contains a positive
 definite matrix}.
\]
Indeed, both conditions are equivalent to the evaluation vectors
$\{v(a)\mid[a]\in Z(h)\}$ spanning $\R^k$.

The following theorem provides a homotopy equivalence between $E$ and $D$. Its determinant maximization is a L\"owner-ellipsoid problem, whose standard contact-point condition is recalled below \cite{John,DanzerLaugwitzLenz}; the SOS/Gram-matrix application is inspired by \cite{PeyrlParrilo}.

\begin{theorem}\label{thm:core}
For $f\in\interior{E}$ there is a unique positive definite matrix $Q_f$ maximizing $\det Q$ subject to
\[
 Q\geqslant 0,\qquad \Gamma(Q)\leqslant f\quad\text{on }S^{n-1}.
\]
The maps $f\mapsto Q_f$ and
\[
 R(f)=\frac{f-\Gamma(Q_f)}{L(f-\Gamma(Q_f))}
\]
are continuous and semialgebraic, and $R(f)\in D$. Both $R\colon\interior{E}\to D$ and the inclusion $D\hookrightarrow E$ are homotopy equivalences. The normalized remainder $R(f)$ is independent of the chosen basis and equivariant under  invertible linear changes of variables, followed by normalization.
\end{theorem}
\begin{proof}
The feasible set is convex and compact by $\tr Q\leqslant L(f)=1$ and contains a positive definite matrix because $f$ is strictly positive. A maximizer is therefore positive definite and unique by strict concavity of $\log\det$.

For continuity, let $f_j\to f$. Closedness of the feasibility conditions implies that every subsequential limit of $Q_{f_j}$ is feasible for $f$. For $0<\eps<1$, the matrix $(1-\eps)Q_f$ is strictly feasible for all sufficiently large $j$. Comparing determinants and then letting $\eps$ go to zero shows that every limit is $Q_f$. Semialgebraicity follows by quantifying the polynomial determinant inequality over the feasible matrices.

For the optimality argument, extend the notation $Q_g$ to arbitrary
strictly positive forms $g$, not necessarily normalized, using the same
determinant-maximization problem. For geometric intuition only, $Q_g$
is the centered L\"owner ellipsoid matrix of  $$K_g=\left\{\pm\frac{v(a)}{\sqrt{g(a)}}\;\middle|\;a\in S^{n-1}\right\}.$$ Indeed, $\Gamma(Q)\leqslant g$
if and only if $K_g\subseteq\{x\in\R^k\mid x^TQx\leqslant1\}$.

We claim that,
for every nonnegative $h$ and $Q>0$,
\begin{equation}\label{eq:kkt}
 Q=Q_{h+\Gamma(Q)}\quad\Leftrightarrow\quad Q^{-1}\in C_h.
\end{equation}
If $Q^{-1}=\sum_{i=1}^s\lambda_i v(a_i)v(a_i)^T\in C_h$, then for every
positive definite $Q'$ feasible for $h+\Gamma(Q)$, the tangent inequality
for $\log\det$ gives
\begin{align*}
 \log\det Q'-\log\det Q
 &\leqslant\tr\bigl(Q^{-1}(Q'-Q)\bigr)\\
 &=\sum_{i=1}^s\lambda_i\Gamma(Q'-Q)(a_i)
 \leqslant\sum_{i=1}^s\lambda_i h(a_i)=0.
\end{align*}
For singular $Q'$ the determinant inequality is automatic, so $Q$ is the
maximizer. Conversely, if $Q=Q_{h+\Gamma(Q)}$ but $Q^{-1}\notin C_h$,
separation from the closed cone $C_h$ gives $D_0\in\Sym_k(\R)$ such that
\[
 \tr(Q^{-1}D_0)>0,
 \qquad \Gamma(D_0)(a)\leqslant0\quad(a\in Z(h)).
\]
Replacing $D_0$ by $D_0-\eta I$ for sufficiently small $\eta>0$ makes the
second inequality strict, while preserving the first. By compactness,
$t\Gamma(D_0)\leqslant h$ on the sphere for all sufficiently small $t>0$:
the inequality holds near $Z(h)$ by strict negativity, while $h$ has a
positive lower bound away from that neighborhood. Choosing $t$ also small
enough that $Q+tD_0>0$ produces a feasible perturbation with
\[
 \left.\frac{d}{dt}\log\det(Q+tD_0)\right|_{t=0}
 =\tr(Q^{-1}D_0)>0,
\]
contradicting optimality and proving \eqref{eq:kkt}.

For $Q=Q_g$, set $h=g-\Gamma(Q)$. Writing the contact vectors as
$x_i=v(a_i)/\sqrt{g(a_i)}$, condition~\eqref{eq:kkt} becomes
\[
 Q^{-1}=\sum_i\mu_i x_ix_i^T,\qquad x_i^TQx_i=1,\qquad \sum_i\mu_i=k,
 \qquad \mu_i>0.
\]
For context, note that this is the contact decomposition for the centered L\"owner ellipsoid: the weights $\mu_i/k$ give
$M=Q^{-1}/k$ and satisfy $\max_{x\in K_g}x^TM^{-1}x\leqslant k$, with equality on their support,
the Kiefer--Wolfowitz criterion for D-optimality \cite{DanzerLaugwitzLenz, John, KieferWolfowitz}.

Put $r=f-\Gamma(Q_f)$. This is nonzero and non-sos, since otherwise $f$
would be sos. Since $f=r+\Gamma(Q_f)$, criterion~\eqref{eq:kkt} gives
$Q_f^{-1}\in C_r$. Thus the evaluation vectors on $Z(r)$ span $\R^k$, and
normalizing $r$ proves that $R(f)\in D$.

For $h\in D$, the only positive semidefinite matrix satisfying $\Gamma(Q)\leqslant h$ is $Q=0$, since it is the only such matrix that annihilates the spanning set $\{v(a)\mid a\in Z(h)\}$. The same compactness argument shows that, for $f_j\to h$ any subsequential limit of $Q_{f_j}$ is feasible for $h$, and thus zero. Consequently, the map $f\mapsto Q_f$ extends continuously to $Y=\interior{E}\cup D$ by setting $Q_h=0$ for $h\in D$. The homotopy
\[
 H(f,t)=\frac{f-t\Gamma(Q_f)}{L(f-t\Gamma(Q_f))},\qquad 0\leqslant t\leqslant1,
\]
is a strong deformation retraction of $Y$ onto $D$. For $f\in\interior{E}$ and $t<1$ its numerator is strictly positive and not sos.

Now use the coordinates of \Cref{prop:collar}. Since
$D\subseteq\partial P\setminus F$, the space $Y$ has the form
\[
 \bigl(\Omega\times(0,1)\bigr)\cup\bigl(D'\times\{1\}\bigr)
 \subseteq\Omega\times(0,1]
\]
for a subset $D'\subseteq\Omega$. The homotopy
\[
 (u,s)\mapsto\bigl(u,(1-t/2)s\bigr),\qquad 0\leqslant t\leqslant1,
\]
pushes both $E$ and $Y$ into $\interior{E}$ and restricts to a homotopy
within each of them. Hence the inclusions
$\interior{E}\hookrightarrow Y\hookrightarrow E$ are homotopy
equivalences. The deformation retraction above is a homotopy equivalence
$Y\to D$, whose restriction to $\interior{E}$ is $R$. Thus
$R\colon\interior{E}\to D$ and the inclusion $D\hookrightarrow E$ are
homotopy equivalences.

A change of basis acts on Gram matrices by congruence and scales every feasible determinant by the same positive constant. Hence the maximizing sos form is independent of the basis. An invertible change of variables also transports the feasible sets by congruence, followed by normalization. Uniqueness then gives the equivariance of $R$.
\end{proof}

The contact cone also describes the fibers of the remainder map $R$. This allows us to prove the dimension bound on the homotopy type of $E$ from the main theorem. 

\begin{theorem}\label{thm:fibers}
For $h\in D,$ put $K_h=C_h\cap\{B>0,\ \tr B=1\}$. Then
\begin{equation}\label{eq:fiber}
 R^{-1}(h)\cong K_h\times(0,1).
\end{equation}
In particular every fiber is contractible, of dimension $\rho(h)$, where
\begin{equation}\label{eq:rho}
 \rho(h)=\dim\bigl(H_{n,2d}/I_{2d}(Z(h))\bigr)\geqslant k.
\end{equation}
For the semialgebraic set $D_r=\{h\in D\mid\rho(h)=r\}$ one has
\[
 \dim D_r\leqslant m-r-1.
\]
Consequently $\dim D\leqslant m-k-1$, and $E$ is homotopy equivalent to a finite CW complex of dimension at most $m-k-1$.
\end{theorem}
\begin{proof}
By criterion~\eqref{eq:kkt}, the positive definite matrices occurring as
optimizers over the remainder $h$ are exactly those with $Q^{-1}\in C_h$.
Write any such matrix uniquely as $Q=tB^{-1}$ with $B\in K_h$ and $t>0$,
and define
\[
 \tau_h(B)=\inf\{t>0\mid h+t\Gamma(B^{-1})\in\Sos_{n,2d}\}.
\]
$\tau_h(B)$ is positive, while finiteness follows
from $\Gamma(B^{-1})\in\interior(\Sos)$. The threshold is continuous since points below it lie in the open complement of the sos cone, while
points strictly above it lie in the interior of that cone. This also proves
joint continuity when $h$ varies in $D$.

Normalize $h+t\Gamma(B^{-1})$ by $L$. Criterion~\eqref{eq:kkt} and uniqueness of the optimizer show that this parametrizes the fiber bijectively and continuously with continuous inverse. Rescaling $t$ by $\tau_h(B)$ proves \eqref{eq:fiber}. The convex set $K_h$ has dimension $\rho(h)-1$.

The span of $C_h$ is the image, under the injective dual Gram map, of the span of degree-$2d$ evaluations on $Z(h)$. This gives the equality in \eqref{eq:rho}. Choose $k$ linearly independent degree-$d$ evaluation vectors; their rank-one squares are linearly independent matrices, proving $\rho(h)\geqslant k$.

The semialgebraic fiber-dimension theorem, which follows from Hardt's
semialgebraic triviality theorem~\cite{Hardt}, now gives
\[
 \dim R^{-1}(D_r)=\dim D_r+r\leqslant\dim\interior{E}=m-1.
\]
Triangulate $(P,D)$ so that $D$ is a union of open simplices. In the barycentric subdivision, $D$ deformation retracts onto the full subcomplex spanned by barycentres $b_\sigma$ of simplices whose relative interiors $\sigma^\circ$ lie in $D$, of dimension at most $\dim D$ \cite{BCR}. Applying \Cref{thm:core} gives the same bound for $E$.
\end{proof}

\section{General lower homology bounds}

We first prove vanishing of high degree cohomology for certain sets of psd matrices constructed from projective semialgebraic sets. Applied to the Veronese embedding, this will prove low degree vanishing of homology of $E$ below.

\begin{proposition}\label{lem:kernel}
Let $V$ be an $N$-dimensional real inner-product space, $N\geqslant2$, and let $X\subseteq\mathbb P(V)$ be compact, semialgebraic, and spanning. If $b=\dim X$ and
\[
 K_X=\{A\geqslant 0\mid\tr A=1,\ \mathbb P(\ker A)\cap X\ne\varnothing\},
\]
then
\[
 H^j(K_X)=0\qquad\text{for }j>(b+1)N-\binom{b+1}{2}-2.
\]
\end{proposition}
\begin{proof}
The proof follows the idea of determinant resolutions from \cite{VassilievResolutions}.
Call a subspace $W$ \emph{$X$-spanned} if $X\cap\mathbb P(W)$ spans $W$.
The
complete flags whose members are all $X$-spanned form a semialgebraic subset
of the usual flag manifold, whose dimension we denote by $\ell(X)$.
We first prove by induction on $N$ that
\begin{equation}\label{eq:flag-bound}
 H^j(K_X)=0\qquad(j>\ell(X)+N-2).
\end{equation}
For $N=2$, the kernel map identifies $K_X$ with $X$, and the space of
complete $X$-spanned flags also equals $X$. Thus
\eqref{eq:flag-bound} follows from the cohomological dimension bound for
semialgebraic sets. For $N>2$, filter $K_X$ by the increasing chain of  compact subsets
\[
 K_i=\{A\in K_X\mid\lambda_{\max}(A)\text{ has multiplicity at least }N-i\},
 \quad 0\leqslant i\leqslant N-1.
\]
Here $K_0=\varnothing$ and $K_{N-1}=K_X$. On $K_i\setminus K_{i-1}$
the maximal-eigenvalue eigenspace $E_{\max}(A)$ has dimension $N-i$.
Let $\operatorname{Gr}_i(V)$ denote the Grassmannian of $i$-dimensional linear subspaces of $V$.
Define
\[
 p_i\colon K_i\setminus K_{i-1}\to\operatorname{Gr}_i(V),
 \qquad p_i(A)=E_{\max}(A)^\perp.
\]
For $U\in\operatorname{Gr}_i(V)$, set $X_U=X\cap\mathbb P(U)$. Relative
to the orthogonal decomposition $V=U\oplus U^\perp$, the fiber
$p_i^{-1}(U)$ consists of
\[
 B\oplus\frac{1-\tr B}{N-i}I_{U^\perp},\quad
 B\geqslant 0,\quad\mathbb P(\ker B)\cap X_U\ne\varnothing,\quad
 \tr B+(N-i)\lambda_{\max}(B)<1.
\]
For $i=1$ the fiber is a point or empty, while for $i\geqslant2$, radial rescaling identifies it with the open cone
$(K_{X_U}\times[0,1))/(K_{X_U}\times\{0\})$. Therefore compactly supported cohomology is the reduced cohomology of the one-point compactification. Here the one-point compactification of the open cone is the suspension of \(K_{X_U}\), so the suspension isomorphism gives
\[
 H_c^j\bigl(p_i^{-1}(U)\bigr)\cong\tH^{j-1}(K_{X_U}).
\]
By Hardt triviality \cite{Hardt} we may partition $\operatorname{Gr}_i(V)$ into finitely many
cells $C_1,\ldots,C_r$ on which both these fibers and the spaces of complete
$X_U$-spanned flags are trivial. By semialgebraic triangulation we can assume that every partial
union $C_1\cup\cdots\cup C_a$ is closed.
If $U$ is not $X$-spanned, then
$X_U$ is either empty or spans a proper subspace of $U$. In the latter
case we can choose a trace-one psd matrix $A_0$ whose kernel contains
$\Span(X_U)$. For any $A\in K_{X_U}$ the straight line path $(A,A_0)$ then clearly is contained in $K_{X_U}$ and so $K_{X_U}$ is contractible. If
$X_U$ is empty, the fiber is empty so in either case the corresponding
open-cone fiber has zero compactly supported cohomology.

Now let $C$ be an $s$-dimensional cell of $X$-spanned subspaces, so $C\cong\R^s$.
Triviality over $C$ identifies $p_i^{-1}(C)$ with $C\times p_i^{-1}(U)$ for any $U\in C$.
Since $H_c^*(\R^s;\Z)$ is concentrated in degree $s$, the compactly supported product formula gives, for $i\geqslant2$,
\[
 H_c^q\bigl(p_i^{-1}(C)\bigr)
 \cong H_c^{q-s}\bigl(p_i^{-1}(U)\bigr)
 \cong \tH^{q-s-1}(K_{X_U}).
\]
The induction hypothesis therefore shows that cohomology can occur only
in degrees
\[
 q\leqslant s+\ell(X_U)+i-1
 \leqslant\ell(X)+N-2.
\]
For the second inequality, we use that the family of such flags over $C$ has dimension
$s+\ell(X_U)$ and every member extends to a complete $X$-spanned flag.
Hence $s+\ell(X_U)\leqslant\ell(X)$, while $i-1\leqslant N-2$. For $i=1$
a nonempty fiber is a point, so over $C$ its
compactly supported cohomology is concentrated in degree $s$, and the same
extension argument gives $s\leqslant\ell(X)$.

It remains to combine the cellwise bounds. Set
\[
 B_a=\bigcup_{b=1}^a C_b,
 \qquad Y_a=p_i^{-1}(B_a).
\]
Then
\[
 \varnothing=Y_0\subseteq Y_1\subseteq\cdots\subseteq Y_r
 =K_i\setminus K_{i-1},
 \qquad Y_a\setminus Y_{a-1}=p_i^{-1}(C_a),
\]
is a closed filtration. At each stage, the compactly supported exact
sequence contains
\[
 \cdots\to H_c^q(Y_a\setminus Y_{a-1})
 \to H_c^q(Y_a)\to H_c^q(Y_{a-1})
 \to\cdots .
\]
The cellwise bounds above and induction on $a$ therefore give
\[
 H_c^q(K_i\setminus K_{i-1})=0
 \qquad(q>\ell(X)+N-2).
\]
Applying the same argument to the closed filtration
$K_i$, whose members are compact and whose successive differences satisfy the preceding bound, finally gives \eqref{eq:flag-bound}: \[H^q(K_X)=H_c^q(K_X)=0 \qquad(q>\ell(X)+N-2).\]
We are only left to bound $\ell(X)$ depending on $b$. The first $N-b-1$ members of a flag can be obtained from $N-b-1$ independent points of $X$; the rest lie in a quotient of dimension $b+1$. Consequently
\[
 \ell(X)\leqslant (N-b-1)b+\frac{b(b+1)}2
 =bN-\binom{b+1}{2}.
\]
This proves the claimed bound.
\end{proof}

This now implies vanishing of homology/cohomology in low/high degree for $E$ and $F$.

\begin{corollary}\label{prop:low-homology}
We have
\[
 H^j(F)=0\quad\left(j>nk-\binom n2-2\right),\qquad
 \tH_i(E)=0\quad(0\leqslant i\leqslant c).
\]
\end{corollary}
\begin{proof}
Apply \Cref{lem:kernel} to the degree-$d$ Veronese image
\[
 X=\{[v(a)]\mid[a]\in\RP^{n-1}\}\subseteq\RP^{k-1},
\]
which spans and has dimension $n-1$. The Gram map sends $K_X$ onto $F,$ since a positive semidefinite Gram matrix represents an sos form vanishing at $a$ exactly when $v(a)$ is in its kernel. Its fibers are nonempty compact convex sets, so the Vietoris--Begle theorem for \v Cech cohomology with integral coefficients \cite{Begle} gives
\[
 \Gamma^*\colon \check H^*(F;\Z)\xrightarrow{\cong}\check H^*(K_X;\Z).
\]
Since the spaces are compact semialgebraic, their \v Cech and singular cohomology groups agree. This gives the first vanishing statement, Alexander duality gives the second.
\end{proof}

Combining this corollary with \Cref{thm:fibers} gives
\[
 \tH_i(E)=0
 \qquad\text{for}\qquad
 0\leqslant i\leqslant c\quad\text{or}\quad i\geqslant m-k.
\]
The lower vanishing range proves connectedness of $E$ whenever $c\geqslant0$.
It remains to examine the cases not covered by this bound. Note that
\[
 \frac{m}{k}=\prod_{j=1}^{n-1}\frac{2d+j}{d+j}
\]
increases with $d$. It is at least $n$ for $n=3,d\geqslant4$;
$n=4,d\geqslant3$; $n=5,d\geqslant3$; and $n\geqslant6,d\geqslant2$.
For the last assertion, at $d=2$ it is
$(n+2)(n+3)/12\geqslant n$. Also $c(5,2)=3$. Thus $c\geqslant0$
outside the two minimal nonempty cases $(n,d)=(3,3),(4,2)$, and the empty cases. \Cref{prop:low-homology}
gives connectedness in all those cases, and local path-connectedness of
semialgebraic sets gives path-connectedness. It remains to treat the two cases
$(n,d)=(3,3)$ and $(4,2)$.

\section{Path-connectedness}
To obtain path connectedness in all cases, we need to study the two minimal nonempty cases separately. Here, our argument is geometric and follows the initial idea of Hilbert: for a
finite zero configuration we compare the forms singular on the configuration
with the span of products of forms vanishing there. The seven-point ternary result from \cite[Theorem~4.3]{ReznickHilbert} has exactly the incidence
hypotheses used below, and the analysis from  \cite{Blekherman}
 gives a complementary explanation of the same two
critical degrees .

For an unordered finite configuration $\mathcal Z \subseteq \mathbb{RP}^{n-1}$, write
\[
 U_{\mathcal Z}=I_d(\mathcal Z),\quad V_{\mathcal Z}=I^{(2)}_{2d}(\mathcal Z),\quad
 W_{\mathcal Z}=\Span(U_{\mathcal Z}U_{\mathcal Z}),\quad
 P_{\mathcal Z}=\{p\in P\mid p|_{\mathcal Z}=0\}.
\]
Here $I^{(2)}$ means vanishing to order at least two. Nonnegativity gives $P_{\mathcal Z}\subseteq V_{\mathcal Z}$, and $S\cap P_{\mathcal Z}\subseteq W_{\mathcal Z}$ by \cite[Theorem~3.2]{ReznickHilbert}. Call $\mathcal Z$ \emph{admissible} if degree-$d$ evaluation is onto, the common real projective zeros of $U_{\mathcal Z}$ are exactly $\mathcal Z$, and its differentials span the projective cotangent space at each point of $\mathcal Z$. The interpolation conditions parallel $d$-independence from \cite{BIK}.
  
The next lemma shows that non-sos forms remain path-connected as their prescribed zeros vary through an admissible family.

\begin{lemma}\label{lem:config-path}
Suppose either $(n,d,r)=(3,3,7)$ or $(n,d,r)=(4,2,6)$. If
$\mathcal A$ is a nonempty path-connected, permutation-invariant family of
ordered admissible $r$-point configurations in one affine chart, then
\[
 \bigcup_{\mathcal Z\in\mathcal A}(P_{\mathcal Z}\setminus S)
\]
is path-connected.
\end{lemma}
\begin{proof}
The space $U_{\mathcal Z}$ has dimension $n$. Multiplication $\Sym^2 U_{\mathcal Z}\to V_{\mathcal Z}$ is injective. Indeed, if a quadratic relation $Q(q_1,\ldots,q_n)=0$ held for a basis of $U_{\mathcal Z}$, its quadratic Taylor term at a prescribed zero would force $Q$ to vanish on the hyperplane $\operatorname{im}Dq$. A quadratic vanishing on a hyperplane factors into linear forms, and the integral-domain property then contradicts linear independence of the $q_i$.

One can prescribe independently, at each point of $\mathcal Z$, the value and all first partial derivatives (in the fixed affine chart) of a degree-$2d$ form. Indeed, take degree-$d$ interpolants $c_i(a_j)=\delta_{ij}$; their squares prescribe the values, while the forms $c_iq$, with $q\in U_{\mathcal Z}$, adjust the first derivatives without changing any values. Admissibility ensures that these derivatives can be adjusted arbitrarily. Thus
\[
 \dim V_{\mathcal Z}=\binom{n+1}{2}+1,\qquad
 \dim W_{\mathcal Z}=\binom{n+1}{2}.
\]
These are the dimension-gap quantities studied by \cite{BIK}.
Fix inner products and orientations on the polynomial spaces. Choose affine coordinates on the fixed chart and dehomogenize the forms. The maps recording values and recording values together with all first partial derivatives then take values in the fixed oriented spaces $\R^r$ and $(\R^n)^r$, respectively, so their kernels $U_{\mathcal Z}$ and $V_{\mathcal Z}$ inherit continuous orientations. Orient $W_{\mathcal Z}$ through $\Sym^2 U_{\mathcal Z}$, and let $b_{\mathcal Z}$ be the positive unit vector in $V_{\mathcal Z}\cap W_{\mathcal Z}^\perp$. A transposition of two elements of $\mathcal Z$ changes the orientations of $U_{\mathcal Z},V_{\mathcal Z},W_{\mathcal Z}$ by $-1,(-1)^n,(-1)^{n-1}$, respectively. Hence
\begin{equation}\label{eq:sign-normal}
 b_{\tau\mathcal Z}=-b_{\mathcal Z}\qquad(\tau\text{ a transposition}).
\end{equation}
Let $s_{\mathcal Z}=\sum q_i^2$ for an orthonormal basis of $U_{\mathcal Z}$. As $\mathcal Z$ varies, $s_{\mathcal Z}$ varies continuously: the spaces $U_{\mathcal Z}$ vary continuously, and this sum is independent of the chosen orthonormal basis. It is also independent of the ordering of $\mathcal Z$ and has exactly the prescribed zeros, all with positive definite projective Hessians. For a fixed configuration, $s_{\mathcal Z}+\eps b_{\mathcal Z}$ is the pointwise perturbation from \cite[Lemma~3.1 and Theorem~3.4]{ReznickHilbert}; for the path argument, the additional issue is choosing its sign continuously and $\eps$ uniformly along paths. On a compact path, Taylor estimates near the zeros and compactness away from them give $|b_{\mathcal Z}|\leqslant C s_{\mathcal Z}$ uniformly. Choose $0<\eps<1/C$. Then $s_{\mathcal Z}+\eps b_{\mathcal Z}$ is nonnegative and is non-sos because it lies outside $W_{\mathcal Z}$.

Take two forms in the union and order their respective configurations so that the corresponding $b_{\mathcal Z}$-components are nonnegative, using \eqref{eq:sign-normal}. Join the ordered configurations in $\mathcal A$. At either endpoint, the segment from the given form to $s_{\mathcal Z}+\eps b_{\mathcal Z}$ has positive $b_{\mathcal Z}$-component except possibly at the original endpoint, which is already non-sos. These segments and the configuration path stay nonnegative and non-sos. Now normalize by $L$ to complete the proof.
\end{proof}
 
We now verify admissibility and path-connectedness for the two configuration types
needed below. \cite[Lemma~2.6]{ReznickHilbert} gives the ternary fullness and real-base-locus
criteria, while the related $d$-independence and dimension-gap framework is developed in \cite{BIK}. Our additional steps are the injectivity of $\Sym^2 U_{\mathcal Z}\to
V_{\mathcal Z}$, the orientation and uniform perturbation, and the
quaternary case with four coplanar points, beyond the general-position result of \cite{BIK}.
 
\begin{lemma}\label{lem:admissible}
The following configurations are admissible, and their ordered spaces in any fixed affine chart are path-connected:
\begin{enumerate}
\item seven points in $\RP^2$, with no four collinear and no conic through all seven, for $d=3$;
\item six points in $\RP^3$, with no three collinear and no five coplanar, for $d=2$.
\end{enumerate}
\end{lemma}
\begin{proof}
We check the three admissibility conditions:
surjectivity of degree-$d$ evaluation, equality of the real projective base
locus with $\mathcal Z$, and spanning of the projective cotangent space by the
differentials at each point. We use the following
observation: on a smooth curve, a simple zero supplies a tangent covector,
while a product of a local equation of the curve with a nonvanishing factor
supplies a normal covector.

\textit{Admissibility of \partlabel{i}.}

\textit{Surjectivity of cubic evaluation:}
Fix $a\in\mathcal Z$. Among the other six points, choose two on a line $L$
avoiding $a$, and let $T$ be the remaining four. If no three points of $T$
are collinear, the two reducible conics obtained by pairing the four points in
two different ways show that the conic pencil through $T$ has base locus
$T$. If instead three points lie on a line $M$ and $b$ is the fourth, then
$I_2(T)=M I_1(b)$; moreover, $a\notin M$ because no four points of
$\mathcal Z$ are collinear. Thus in either case there is
$C\in I_2(T)$ with $C(a)\ne0$. The cubic $LC$ vanishes on
$\mathcal Z\setminus\{a\}$ but not at $a$. Thus cubic evaluation is onto and
$\dim U_{\mathcal Z}=10-7=3$.

\textit{Base locus and differential spanning when three points are
collinear:}
Suppose that three points lie on a line $L$, and set
$T=\mathcal Z\setminus L$. Then
\[
 U_{\mathcal Z}=L I_2(T)\oplus\R f,\qquad L\nmid f.
\]
If no three points of $T$ are collinear, the two reducible conics obtained by
pairing the four points in two different ways show that the pencil $I_2(T)$
has base locus $T$ and spans the cotangent spaces there. Moreover, $f|_L$ is
a nonzero binary cubic whose three prescribed zeros are therefore simple and
exhaustive. This proves both the base-locus and differential assertions.

If instead three points of $T$ lie on a line $M$, with fourth point $b$, then
$I_2(T)=M I_1(b)$, so $L I_2(T)=LM I_1(b)$. A cubic in
$U_{\mathcal Z}$ divisible by $M$
must, after division by
$M$, contain $L$, so $f$ is divisible by neither $L$ nor $M$. Its
restrictions to both lines have precisely the three prescribed simple zeros.
Since the base locus of $LM I_1(b)$ is $L\cup M\cup\{b\}$, the base locus of
$U_{\mathcal Z}$ is exactly $\mathcal Z$. The same restrictions give the
tangent covectors on $L$ and $M$. $LM I_1(b)$ gives the normal covectors
there and spans the cotangent space at $b$.

\textit{Base locus and differential spanning when no three points are collinear:}
Any five points determine a unique conic $C$, which is smooth, since a
singular real conic through five such points would contain three of them on a
line. By hypothesis $C$ contains exactly five or six points of $\mathcal Z$.

If $C$ contains five points, let $L$ join the other two. Restriction of
$U_{\mathcal Z}$ to either $C$ or $L$ has kernel $\R CL$. Its image is the
full two-dimensional
space of, respectively, binary sextics vanishing at the five prescribed points of
$C$ or binary cubics vanishing at the two prescribed points of $L$. These systems
have only the prescribed base points and contain sections with a simple zero at
each of them; $CL$ supplies the normal covectors.

If $C$ contains six points and $b$ is the seventh, then
\[
 U_{\mathcal Z}=C I_1(b)\oplus\R f,\qquad C\nmid f.
\]
The restriction $f|_C$ has precisely the six prescribed simple zeros, while
$C I_1(b)$ supplies the normal covectors on $C$ and all cotangent directions
at $b$. Its base locus is $C\cup\{b\}$, so intersecting with $Z(f)$ leaves
exactly $\mathcal Z$. This completes the base-locus and differential checks,
and hence the admissibility of \partlabel{i}.

\textit{Admissibility of \partlabel{ii}.}

\textit{Surjectivity, base locus, and differential spanning when no four points are coplanar:}
Split the six points into two triples spanning planes $\{L=0\}$ and
$\{M=0\}$. Restriction of $U_{\mathcal Z}$ to either plane is onto the
three-dimensional space of conics through the corresponding triple: an
arbitrary lift can be corrected at the other triple by adding $L\ell$ or
$M\ell$. The kernel is $\R LM$, hence
$\dim U_{\mathcal Z}=4=10-6$, which proves that degree-two evaluation is
onto.

On a plane, the conics through three noncollinear points have exactly those
base points and span their cotangent spaces; in suitable coordinates one may
use $xy,xz,yz$. Since $LM\in U_{\mathcal Z}$, every common zero lies on the
union of the two planes. The restricted systems give precisely the six
prescribed zeros and their tangent covectors, while $LM$ gives the normal
covectors. Thus the base-locus and differential conditions also hold.

\textit{Surjectivity, base locus, and differential spanning when four points are coplanar:}
Let $\Pi=\{L=0\}$ contain the four points, let $p,q$ be the other two, and set
$M=\overline{pq}$. Restriction to $\Pi$ is onto the conic pencil through its
four points and has kernel
\[
 L I_1(M)=\{L\ell\mid \ell|_M=0\}.
\]
Both spaces have dimension two, so $\dim U_{\mathcal Z}=4=10-6$ and
degree-two evaluation is onto.

The pencil has precisely the four prescribed base points and spans their
cotangent spaces in $\Pi$. Outside $\Pi$, the kernel forces a common zero to
lie on $M$. Since $M\cap\Pi$ is not one of the prescribed points, choose a conic in the pencil
nonzero there and lift it to $Q\in U_{\mathcal Z}$. Then $Q|_M$ is a nonzero quadratic
with the two distinct zeros $p,q$, hence no others, and both zeros are simple.
The pencil together with $L I_1(M)$ spans the cotangent spaces at the four
plane points, while at $p,q$ the latter supplies the two directions normal to
$M$ and $Q$ supplies the tangent direction. This verifies the base-locus and
differential conditions, so the configuration is admissible.

\textit{Path-connectedness of the ordered configuration spaces.}

For \partlabel{i}, the ambient space has dimension $14$. The locus where two
of the seven points coincide has codimension two; the locus with a collinear quadruple has dimension at most
$2+4+6=12$; and the incidence of a projective conic with seven points on it
has dimension at most $5+7=12$. For \partlabel{ii},
the ambient space has dimension $18$. The locus where two of the six points
coincide has codimension three,
while the loci with a collinear triple or five coplanar points have dimensions
at most
\[
 4+3+9=16,\qquad 3+10+3=16,
\]
respectively. Hence in either case the forbidden locus is semialgebraic and
has codimension at least two, so its complement is path-connected.
\end{proof}

We can now finally prove path connectedness for all cases.
\begin{theorem}\label{thm:connected}
$E$ is path-connected.
\end{theorem}
\begin{proof}
By \Cref{thm:core}, in each of the two remaining cases every form in
$E$ can be joined to some $h\in D$. We construct an admissible configuration
among the zeros of $h$.

By \cite[Proposition~7]{BHORS}, a non-sos ternary sextic or quaternary quartic
cannot vanish on a line. A nonzero nonnegative binary form of degree $2d$ has
at most $d$ distinct real projective zeros, since each has even multiplicity.
Thus four collinear zeros of a ternary sextic, or three of a quaternary
quartic, would force the restriction to the line to vanish identically and
are therefore impossible.

For ternary sextics, $I_3(Z(h))=0$ implies $I_2(Z(h))=0$ and
$|Z(h)|\geqslant\dim H_{3,3}=10$. The quadratic evaluations therefore span
$H_{3,2}^*$; choose six forming a basis and add a seventh zero, as in \cite[Lemma~7.5 and Theorem~7.6]{ReznickHilbert}. No conic
contains the first six, and no four of the seven points are collinear, so
\Cref{lem:admissible}\partlabel{i} applies.

For quaternary quartics, $I_2(Z(h))=0$ implies that $Z(h)$ spans $\RP^3$ and
has at least $\dim H_{4,2}=10$ points. Choose four noncoplanar zeros and a
fifth zero, forming a set $T$. If no four points of $T$ are coplanar, add any
further zero. Otherwise a unique plane $\Pi$ contains four points of $T$;
let $b$ be the remaining point and choose a zero outside $\Pi\cup\{b\}$.
If no such zero existed, then the product of
an equation for $\Pi$ and an equation of any plane through $b$ would belong
to $I_2(Z(h))$, a contradiction. In the
first case a plane through five selected points would contain four points of
$T$; in the second it would contain three noncollinear points of $\Pi$ and a
point outside $\Pi$. Hence no five selected points are coplanar, and
\Cref{lem:admissible}\partlabel{ii} applies.

Given two members of $D$, choose these configurations and a common affine
chart, obtained by choosing a hyperplane avoiding their finite union.
\Cref{lem:admissible} connects the ordered configurations, and
\Cref{lem:config-path} lifts this to a path in $E$.
\end{proof}

\section{Euler characteristic}
For a nonempty space $X$ whose rational homology groups are finite-dimensional
and vanish in sufficiently high degrees, write $\tchi(X)=\chi(X)-1$.
Alexander duality gives
\begin{equation}\label{eq:euler-dual}
 \tchi(E)=(-1)^{m-3}\tchi(F).
\end{equation}
We will compute Euler characteristics by passing to fixed-point sets of
symmetries. If a finite-order transformation $g$ preserves a compact 
semialgebraic set $Y$ in a linear representation, equivariant semialgebraic
triangulation gives $Y$ a finite $\langle g\rangle$-CW structure, so Brown's
periodic Lefschetz formula applies \cite{Brown, ParkSuh}. It identifies the
Euler characteristic of the fixed-point set with the alternating trace of $g$ on the cohomology of $Y$:
\begin{equation}\label{eq:lefschetz}
 \chi(Y^g)=\sum_i(-1)^i\tr\bigl(g^*\colon H^i(Y;\Q)\to H^i(Y;\Q)\bigr).
\end{equation}
In particular, if $g$ belongs to a path-connected group preserving $Y$, its action
on $Y$ is homotopic to the identity, and the right-hand side is $\chi(Y)$.
For a torus $T$ acting linearly and preserving $Y$, complexify the ambient
representation and decompose it into weight spaces. Since only finitely many
weights occur, we can choose a finite-order element $g\in T$ on which every
nontrivial weight takes a value other than $1$. A vector is fixed by $g$
precisely when it lies in the weight-zero subspace, which is also the subspace
fixed by all of $T$. Consequently, $Y^g=Y^T$ and
\begin{equation}\label{eq:torus-euler}
 \chi(Y)=\chi(Y^T).
\end{equation}
This standard torus fixed-point formula \cite{Bredon} lets us compute the Euler
characteristic from an easier fixed point space.
We used very similar  methods in  \cite{IllmerNetzerWolf}  to study  the topology of the set of entangled states. 

\begin{theorem}\label{thm:euler}
With $\chi(\varnothing)=0$,
\begin{enumerate}
\item for even $n\geqslant4$, $\chi(F_{n,2d})=\chi(F_{n-1,2d})$ and $\chi(E_{n,2d})=\chi(E_{n-1,2d})$;
\item $\chi(F_{n,2d})$ and $\chi(E_{n,2d})$ are even for all $n\geqslant2$;
\item $\chi(E_{3,2d})=\chi(E_{4,2d})=0$, and $\chi(E_{n,4})=0$ for $5\leqslant n\leqslant8$.
\end{enumerate}
\end{theorem}
\begin{proof}
Let a circle rotate $(x_1,x_2)$. Invariant polynomials are of the form
\[
 p(x)=h(x_1^2+x_2^2,x_3,\ldots,x_n).
\]
Put $x'=(x_3,\ldots,x_n)$ and rename $x_1$ as $y$. Restriction to
$x_2=0$ sends $p$ to
\[
 \widehat p(y,x')=h(y^2,x'),
\]
which is fixed by the reflection $\tau(y,x')=(-y,x')$. Conversely, every
$\tau$-invariant form contains only even powers of $y$, and hence has a
unique expression of this form; replacing $y^2$ by $x_1^2+x_2^2$ recovers
$p$. Thus restriction is a linear bijection
\[
 H_{n,2d}^{S^1}\to H_{n-1,2d}^{\tau}
\]
from the circle-invariant forms in $n$ variables to the
$\tau$-invariant forms in $n-1$ variables.

The bijection preserves nonnegativity: for every $(x_1,x_2,x')$, putting
$y=(x_1^2+x_2^2)^{1/2}$ gives $p(x_1,x_2,x')=\widehat p(y,x')$. It also
preserves the existence of a nonzero real zero, by the same substitution in
one direction and by evaluating $p(y,0,x')$ in the other.

Restriction preserves sums of squares. Conversely, if
$\widehat p=\sum_i q_i^2$, average this expression under $\tau$ and decompose
each $q_i$ into its even and odd parts in $y$. The even part has the form
$a_i(y^2,x')$, while the odd part has the form $y b_i(y^2,x')$, giving
\[
 \widehat p=\sum_i a_i(y^2,x')^2+y^2\sum_j b_j(y^2,x')^2,
\]
up to relabeling the indices. After substituting
$y^2=x_1^2+x_2^2$, the second term becomes
\[
 \sum_j\bigl(x_1b_j(x_1^2+x_2^2,x')\bigr)^2
 +\sum_j\bigl(x_2b_j(x_1^2+x_2^2,x')\bigr)^2,
\]
so the recovered form $p$ is again sos. After normalization, therefore,
\begin{equation}\label{eq:fixed-reduction}
 F_{n,2d}^{S^1}\cong F_{n-1,2d}^{\tau}.
\end{equation}

If $n$ is even, $-\tau\in\SO(n-1)$ and induces the same action as $\tau$
on even-degree forms. Since $\SO(n-1)$ is path-connected, this action is
homotopic to the identity. Applying \eqref{eq:lefschetz} to $\tau$ on
$F_{n-1,2d}$ gives $\chi(F_{n-1,2d}^{\tau})=\chi(F_{n-1,2d})$. On the
other hand, applying \eqref{eq:torus-euler} to $F_{n,2d}$ and then using
\eqref{eq:fixed-reduction} gives
$\chi(F_{n,2d})=\chi(F_{n-1,2d}^{\tau})$. Thus
$\chi(F_{n,2d})=\chi(F_{n-1,2d})$. Moreover
\[
 m_n-m_{n-1}=\binom{2d+n-2}{2d-1}
\]
is even, so \eqref{eq:euler-dual} proves the equality for $E$ whenever both
spaces are nonempty. If $d=1$, both $E_{n,2}$ and $E_{n-1,2}$ are empty, so
the equality is trivial. For $(n,2d)=(4,4)$,
$\Pos_{3,4}=\Sos_{3,4}$ by Hilbert's theorem, and hence
$F_{3,4}=\partial P_{3,4}$. Since
$\dim P_{3,4}=\dim H_{3,4}-1=15-1=14$, the boundary of this compact convex
body is homeomorphic to $S^{13}$. Thus $\chi(F_{3,4})=0$, and the same formulas give
$\chi(E_{4,4})=0=\chi(E_{3,4})$.

For odd $n$, the eigenvalues of the involution $\tau^*$ are $\pm1$, so its trace is congruent modulo two to the dimension in every degree. Thus \eqref{eq:fixed-reduction} implies
$\chi(F_{n,2d})\equiv\chi(F_{n-1,2d})\pmod2$. The base case is $F_{2,2d}\cong S^{2d-1}$, of Euler characteristic zero. This proves evenness for $F$, and \eqref{eq:euler-dual} proves it for every $E$.

For $n=3$, circle-invariant forms have the form
\[
 \sum_{j=0}^{d}a_j x_1^{2d-2j}(x_2^2+x_3^2)^j.
\]
Every nonnegative such form is sos, by restricting to an even binary form and applying the preceding averaging argument. The normalized invariant cone is a convex body of dimension $d$, whose boundary is exactly $F^{S^1}$. Hence $F^{S^1}\cong S^{d-1}$ and $\chi(F)=1+(-1)^{d-1}$. For $d\leqslant2$ the ternary space $E$ is empty. For $d\geqslant3$, since $m+d$ is odd, \eqref{eq:euler-dual} gives $\chi(E)=0$. Part \partlabel{i} gives the quaternary assertion.

For quartics in $n=6$ or $8$ variables, let $T=(S^1)^r$, where $r=n/2$, act by independently rotating the coordinate pairs. Invariant quartics are of the form $u^TAu$, where $u_i=x_{2i-1}^2+x_{2i}^2$ and $A$ is copositive. For $r\leqslant4$, \cite{Diananda} writes $A$ as a positive semidefinite matrix plus an entrywise nonnegative matrix. Thus every invariant nonnegative quartic is sos. The normalized cone of invariant nonnegative quartics is a convex body of dimension $r(r+1)/2-1$, and its boundary is exactly $F^T$. Hence $F^T\cong S^{r(r+1)/2-2}$, of dimension $4$ or $8$. The fixed-point Euler-characteristic formula gives $\chi(F)=\chi(F^T)=2$. Since $m=126$ or $330$ is even, \eqref{eq:euler-dual} gives $\chi(E)=0$. Part \partlabel{i} gives the cases $n=5,7$.
\end{proof}

$E$ is connected, but its even Euler characteristic rules out rational acyclicity, i.e.\ the vanishing of all rational reduced homology groups. When $\chi(E)=0$, the alternating sum of positive-degree Betti numbers is $-1$, so at least one odd-degree group has positive rank. Combined with the results of Sections 2--4, this proves \Cref{thm:main-general}.

\section{The two minimal nonempty cases}\label{sec:contact-covers}

We now treat ternary sextics and quaternary quartics in more detail. In this section
\[
 X_{3}= \interior{E_{3,6}},\qquad X_{4}=\interior{E_{4,4}}.
\]
The inclusions \(X_3\hookrightarrow E_{3,6}\) and \(X_4\hookrightarrow E_{4,4}\) are homotopy equivalences by
\Cref{prop:collar}.
For a topological space \(M\) and an integer \(r\geqslant1\), define the ordered configuration space and its unordered quotient by
\[
 \Conf_r(M)=\{(x_1,\ldots,x_r)\in M^r\mid x_i\ne x_j\text{ for }i\ne j\},
 \qquad \UConf_r(M)=\Conf_r(M)/\mathfrak S_r.
\]
We give \(\Conf_r(M)\) the subspace topology and \(\UConf_r(M)\) the quotient topology, with \(\mathfrak S_r\) acting by coordinate permutations. For \(1\leqslant r\leqslant n\), write
\[
 \Conf^{\mathrm{ind}}_r(\RP^{n-1})=
 \{([a_1],\ldots,[a_r])\in\Conf_r(\RP^{n-1})\mid
 a_1,\ldots,a_r\text{ are linearly independent in }\R^n\}.
\]
This condition is independent of the chosen representatives, and it defines an open subspace. 

The upper bound in the ten-zero assertion below is known from \cite{BHORS}: every non-sos nonnegative ternary sextic or quaternary quartic has
at most ten distinct real projective zeros.  The theorem shows that the residual forms attain this bound and that their zero
configurations vary continuously.

\begin{theorem}\label{thm:contact-covers}
For \(X=X_3\) or \(X=X_4\), every \(R(f)\), \(f\in X\), has exactly ten
distinct real projective zeros, and
\[
 \zeta\colon X\to\UConf_{10}(\RP^{n-1}),
 \qquad f\mapsto Z(R(f)),
\]
is continuous.  Consequently, for \(1\leqslant r\leqslant10\),
\[
 \mathcal C_r(X)=
 \{(f,a_1,\ldots,a_r)\mid a_i\in Z(R(f)),\ a_i\ne a_j\}
\]
is a finite covering
\[
 p_r\colon\mathcal C_r(X)\to X
\]
of degree \(N_r=10!/(10-r)!\), and the projection to the ordered root configuration is continuous
\[
 \eta_r\colon\mathcal C_r(X)\to
 \Conf_r(\RP^{n-1}),\qquad
 (f,a_1,\ldots,a_r)\mapsto(a_1,\ldots,a_r).
\]
For \(X=X_4\), the image of \(\eta_3\) lies in the open subspace
\(\Conf^{\mathrm{ind}}_3(\RP^3)\) of linearly independent triples.
\end{theorem}
\begin{proof}
For both $E_{3,6}$ and $E_{4,4}$ we have $k=\dim H_{n,d}=\binom{n+d-1}{d}=10$.  By definition of $D$ we have \(I_d(Z(R(f)))=0\) and so $R(f)$ must have at least $k=10$ distinct zeros.  On the
other hand, a non-sos nonnegative ternary sextic or quaternary quartic has at
most ten distinct real projective zeros \cite[Proposition~7]{BHORS}.  Thus it
has exactly ten.

Let \(f_j\to f\) in \(X\).  Criterion~\eqref{eq:kkt} gives
\[
 Q_{f_j}^{-1}=\sum_{i=1}^{10}\lambda_{ji}
 v(a_{ji})v(a_{ji})^T,\qquad \lambda_{ji}>0,
\]
where the \(a_{ji}\) represent the ten zeros of \(R(f_j)\).  Indeed, the
ten evaluation vectors form a basis, and omitting one generator could not
produce a positive definite matrix.  The matrices \(Q_{f_j}^{-1}\) converge
to \(Q_f^{-1}\), and the coefficients are bounded by taking traces.  After
passing to a subsequence, all points and coefficients converge.  If two
of these projective zeros collided, the limiting positive definite matrix would be a
sum of at most nine rank-one matrices, which is impossible.  The ten limit
points are therefore distinct zeros of \(R(f)\), hence all of its zeros.
This proves continuity of \(\zeta\).  The asserted coverings are pullbacks along $\zeta$
of the standard finite coverings that choose an ordered \(r\)-tuple from an
unordered ten-point configuration.

It remains to treat triples in the quartic case.  Three collinear zeros
would force the restriction of \(R(f)\) to that line to vanish identically,
since a nonnegative binary quartic cannot have three distinct projective
zeros.  A nonnegative quaternary quartic vanishing on a line has infinitely many
zeros, so it is sos by \cite[Proposition~7]{BHORS}. This contradicts
\(R(f)\in D\subseteq E\).
Thus every triple of distinct points of \(Z(R(f))\) is linearly independent.
\end{proof}

The two configuration spaces used in the calculations below have
particularly simple homotopy types.  The orthogonalization idea in the first
part of the lemma is closely related to the deformation
retraction of independent spherical configurations onto a Stiefel manifold
from \cite{WangZhao}.
 
\begin{lemma}\label{lem:small-configurations}
There are deformation retractions
\[
 \Conf_2(\RP^2)\simeq \SO(3)/K_2,\qquad
 \Conf^{\mathrm{ind}}_3(\RP^3)\simeq \SO(4)/K_3,
\]
where \(K_2\) and \(K_3\) are the full stabilizers in \(\SO(3)\) and
\(\SO(4)\), respectively, of an ordered orthogonal projective frame.  For
an orthogonal \(r\)-configuration \(\mathbf a_0\) (with $r=2,3$), the orbit map
\[
 q_{\mathbf a_0}\colon\SO(n)\to
 \Conf^{\mathrm{ind}}_r(\RP^{n-1}),
 \qquad g\mapsto g\mathbf a_0,
\]
induces the isomorphisms
\[
 \begin{aligned}
 q_{\mathbf a_0}^*\colon
 H^*(\Conf_2(\RP^2);\Q)&\xrightarrow{\ \cong\ }
 H^*(\SO(3);\Q)\cong\Lambda_{\Q}(z_3),\\
 q_{\mathbf a_0}^*\colon
 H^*(\Conf^{\mathrm{ind}}_3(\RP^3);\Q)&\xrightarrow{\ \cong\ }
 H^*(\SO(4);\Q)\cong\Lambda_{\Q}(z_3,w_3).
 \end{aligned}
\]
\end{lemma}
\begin{proof}
Choose unit representatives of the ordered projective points and put them
in the columns of a matrix \(A\).  In either case \(A\) has full column rank
and thus \(A^TA\) is positive definite.  The homotopy
\[
 A\mapsto A(A^TA)^{-t/2},\qquad 0\leqslant t\leqslant1,
\]
is unchanged projectively by changing the signs of the columns, is
rotation-equivariant, and retracts the space of configurations onto the subspace of
orthogonal ones.  These are \(\SO(3)/K_2\) and
\(\SO(4)/K_3\), respectively, showing the first statement.

Fix an orthogonal \(r\)-configuration \(\mathbf a_0\).  Its orbit map onto the
subspace of orthogonal \(r\)-configurations is the corresponding finite
covering \(\SO(3)\to\SO(3)/K_2\) or
\(\SO(4)\to\SO(4)/K_3\).  Since \(\SO(n)\) is path-connected, right translation by each
element of \(K_r\) is homotopic to the identity.  Hence \(K_r\) acts
trivially on \(H^*(\SO(n);\Q)\), and
\cite[Proposition~3G.1]{Hatcher} implies that this covering induces an
isomorphism in rational cohomology.  The map \(q_{\mathbf a_0}\) with target
the full configuration space is the composite of this covering with the
inclusion of the orthogonal configurations, so the deformation retraction shows
that it too induces an isomorphism.  Together with the standard
rational cohomology calculations for \(\SO(3)\) and \(\SO(4)\)
\cite[Section~3.D]{Hatcher}, this proves the displayed cohomology formulas.
\end{proof}

Going back to the setup of \Cref{thm:contact-covers}, let
\(\operatorname{Tr}_{p_r}\) denote the cohomological transfer for the finite
cover \(p_r\), as in \cite[Section~3.G]{Hatcher}.  For
\((n,X,r)=(3,X_3,2)\) or \((4,X_4,3)\), define
\[
 T\colon H^*(\Conf^{\mathrm{ind}}_r(\RP^{n-1});\Q)
 \to H^*(X;\Q),\qquad
 \gamma\mapsto
 \frac1{N_r}\operatorname{Tr}_{p_r}(\eta_r^*\gamma).
\]

\begin{proposition}\label{prop:transfer-splitting}\label{prop:orbit-pi1}
For every \(f\in X\), the rotation orbit map
\[
 o_f\colon\SO(n)\to X,
 \qquad g\mapsto g\cdot f=\bigl(x\mapsto f(g^{-1}x)\bigr),
\]
induces a surjection of rational cohomology algebras that admits a
graded-algebra section.
\\
Moreover, it also induces an injection
\[
 (o_f)_*\colon\pi_1(\SO(n))\cong\Z/2\hookrightarrow\pi_1(X,f)
\]
whose image is central.

In the cases $(n,X,r)=(3,X_3,2)$ and $(n,X,r)=(4,X_4,3)$, for every $j\geqslant2$, the induced homomorphism
\[
 (o_f)_*\colon\pi_j(\SO(n),e)\longrightarrow\pi_j(X,f)
\]
is split injective.
\end{proposition}
\begin{proof}
Equivariance of \(R\) gives \(Z(R(g\cdot f))=gZ(R(f))\).  Put
\[
 \mathcal B_f\coloneqq
 \{\text{ordered \(r\)-tuples of distinct points of \(Z(R(f))\)}\},
\]
so \(\lvert\mathcal B_f\rvert=N_r\).  The pullback of \(p_r\) along
\(o_f\) and all the maps used below therefore fit into the diagram
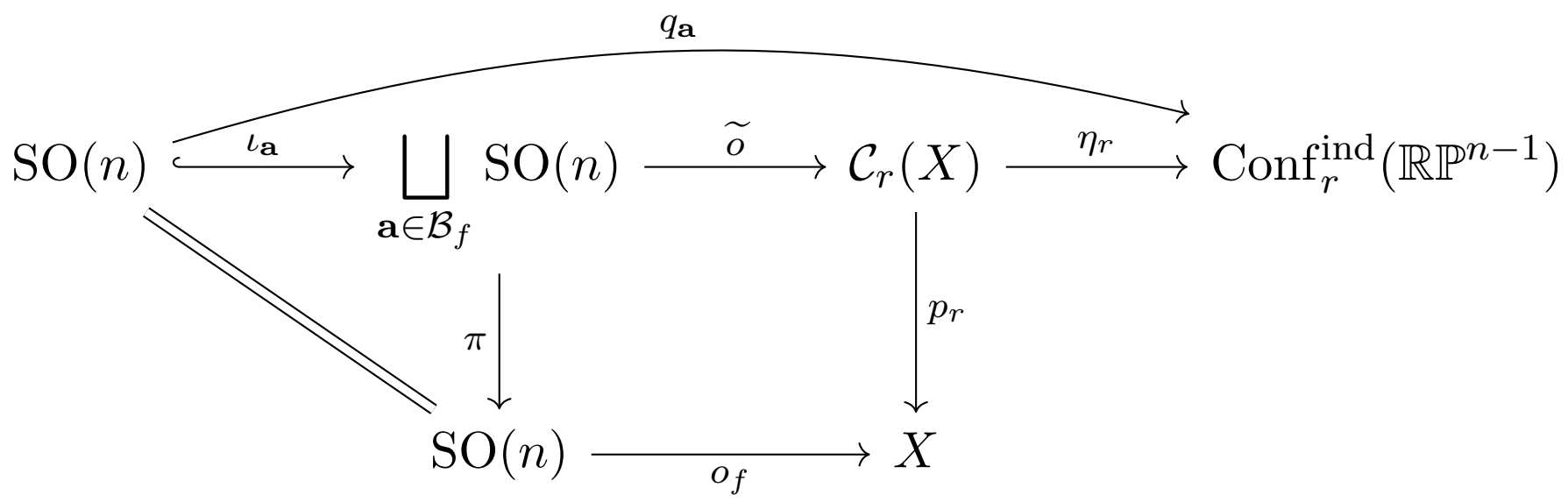
\[
\begin{tikzcd}[column sep=large,row sep=large]
 \SO(n)
   \arrow[r,hook,"\iota_{\mathbf a}"]
   \arrow[dr,equal]
   \arrow[rrr,bend left=15,"q_{\mathbf a}"]
 & \displaystyle\bigsqcup_{\mathbf a\in\mathcal B_f}\SO(n)
   \arrow[r,"\widetilde o"]
   \arrow[d,"\pi"']
 & \mathcal C_r(X) \arrow[r,"\eta_r"] \arrow[d,"p_r"]
 & \Conf^{\mathrm{ind}}_r(\RP^{n-1}) \\
 & \SO(n) \arrow[r,"o_f"']
 & X
 & {}
\end{tikzcd}
\]
Here \(\pi(g,\mathbf a)=g\), \(\iota_{\mathbf a}(g)=(g,\mathbf a)\), and
\(\widetilde o(g,\mathbf a)=(g\cdot f,g\mathbf a)\).  Thus the restriction
\(\widetilde o_{\mathbf a}=\widetilde o\circ\iota_{\mathbf a}\) sends \(g\)
to \((g\cdot f,g\mathbf a)\), while
\(q_{\mathbf a}=\eta_r\circ\widetilde o_{\mathbf a}\) sends \(g\) to
\(g\mathbf a\).
The square with vertical maps \(\pi\) and \(p_r\) is a pullback square,
expressed as the disjoint union of its \(N_r\) sheets.  Naturality of transfer
\cite[Section~3.G]{Hatcher} gives
\[
 o_f^*\circ T
 =\frac1{N_r}\sum_{\mathbf a\in\mathcal B_f}
   \widetilde o_{\mathbf a}^*\circ \eta_r^*
 =\frac1{N_r}\sum_{\mathbf a\in\mathcal B_f}q_{\mathbf a}^*.
\]
Pick an arbitrary orthogonal \(r\)-configuration \(\mathbf a_0\).  The
configuration space is path-connected by
\Cref{lem:small-configurations}.  Thus, for each \(\mathbf a\in\mathcal B_f\),
a path \(\mathbf a(t)\) from \(\mathbf a\) to \(\mathbf a_0\) gives a
homotopy \((g,t)\mapsto g\mathbf a(t)\) from \(q_{\mathbf a}\) to
\(q_{\mathbf a_0}\).  Consequently,
\[
 o_f^*\circ T=q_{\mathbf a_0}^*
\]
where by \Cref{lem:small-configurations} $q_{\mathbf a_0}^*$ is an isomorphism.
Although \(T\) need not be
multiplicative, we can use it to lift the algebra generators.  Write
\[
 H^*(\SO(n);\Q)=\Lambda_{\Q}(z_1,\ldots,z_s),
 \qquad |z_i|=3,\qquad s=1\ (n=3),\;2\ (n=4),
\]
and set
\[
 \widehat z_i=T\bigl((q_{\mathbf a_0}^*)^{-1}(z_i)\bigr)\in H^3(X;\Q).
\]
Then \(o_f^*(\widehat z_i)=z_i\).  Since the \(\widehat z_i\) have odd
degree, graded commutativity gives precisely the exterior-algebra relations.
Thus \(z_i\mapsto\widehat z_i\) extends to a graded-algebra homomorphism
\(s_f\colon H^*(\SO(n);\Q)\to H^*(X;\Q)\), and \(o_f^*s_f\) fixes every
generator.  Hence \(o_f^*\circ s_f=\operatorname{id}\).

For the fundamental-group assertion, choose an ordered \(r\)-tuple
\(\mathbf a\in\mathcal B_f\).  Again \(q_{\mathbf a}\) is
homotopic to \(q_{\mathbf a_0}\) for an orthogonal \(r\)-configuration
\(\mathbf a_0\) and \Cref{lem:small-configurations} shows that the latter
is a finite covering followed by a homotopy equivalence, so
\((q_{\mathbf a})_*\) is injective.  Since \((p_r)_*\) is also injective, applying $\pi_1$ to the above diagram
shows that \((o_f)_*\) is injective.

For any action of a connected group, the image is central: if \(g_s\) is a
loop in the group and \(x_t\) is a loop in the space with basepoint $f$, the map
\((s,t)\mapsto g_s\cdot x_t\) gives a null-homotopy of the commutator of
their images.

For the higher homotopy groups, fix $\mathbf a\in\mathcal B_f$, and use the basepoints $e$, $(f,\mathbf a)$, $f$, and $\mathbf a$ in $\SO(n)$, $\mathcal C_r(X)$, $X$, and $\Conf^{\mathrm{ind}}_r(\RP^{n-1})$, respectively. If $\mathcal C_r(X)$ is disconnected, we take its based homotopy groups in the component containing $(f,\mathbf a)$. Since $p_r$ is a covering, $(p_r)_*$ is an isomorphism for $j\geqslant2$ \cite[Proposition~4.1]{Hatcher}.

Let $\rho_t$ be the equivariant deformation retraction of \Cref{lem:small-configurations} onto the orthogonal-configuration subspace, and put $\mathbf a_0=\rho_1(\mathbf a)$. The homotopy $(g,t)\mapsto g\rho_t(\mathbf a)$ carries $q_{\mathbf a}$ to the orbit map of $\mathbf a_0$, with basepoint path $t\mapsto\rho_t(\mathbf a)$. The latter orbit map onto the orthogonal-configuration subspace is the finite covering $\SO(n)\to\SO(n)/K_r$, which induces an isomorphism on $\pi_j$ for $j\geqslant2$ by \cite[Proposition~4.1]{Hatcher}. The inclusion of that subspace into $\Conf^{\mathrm{ind}}_r(\RP^{n-1})$ is a homotopy equivalence by \Cref{lem:small-configurations}. Thus, accounting for the indicated basepoint path, $(q_{\mathbf a})_*$ is an isomorphism for $j\geqslant2$.

Define
\[
\ell_j=(q_{\mathbf a})_*^{-1}\circ(\eta_r)_*\circ(p_r)_*^{-1}\colon
\pi_j(X,f)\longrightarrow\pi_j(\SO(n),e).
\]
The identities $p_r\circ\widetilde o_{\mathbf a}=o_f$ and $\eta_r\circ\widetilde o_{\mathbf a}=q_{\mathbf a}$ give
\[
\ell_j\circ(o_f)_*=(q_{\mathbf a})_*^{-1}\circ(\eta_r)_*\circ(\widetilde o_{\mathbf a})_*=\operatorname{id}.
\]
Hence $(o_f)_*$ is split injective on the integral homotopy groups.
\end{proof}

\begin{corollary}\label{cor:small-orbit-consequences}%
\label{cor:ternary-pi3}\label{cor:quartic-pi3}
\begin{enumerate}
\item For ternary sextics there is an embedding
\[
 \Lambda_{\Q}(u_3)\hookrightarrow H^*(E_{3,6};\Q),
\]
which implies \(b_3(E_{3,6})\geqslant1\).  In addition, the rotation
orbit map precomposed with the spin covering \(\operatorname{Spin}(3)\to\SO(3)\) induces a split injection
\[
 \Z=\pi_3(\operatorname{Spin}(3))\hookrightarrow\pi_3(E_{3,6}),
\]
and the center of \(\pi_1(E_{3,6})\) contains \(\Z/2\).
\item For quaternary quartics there is an embedding
\[
 \Lambda_{\Q}(u_3,v_3)\hookrightarrow H^*(E_{4,4};\Q),
\]
which implies \(b_3(E_{4,4})\geqslant2\) and
\(b_6(E_{4,4})\geqslant1\).  In addition, the rotation orbit map precomposed with the spin covering
\(\operatorname{Spin}(4)\to\SO(4)\) induces a split injection
\[
 \Z^2=\pi_3(\operatorname{Spin}(4))\hookrightarrow\pi_3(E_{4,4}),
\]
and the center of \(\pi_1(E_{4,4})\) contains \(\Z/2\).
\end{enumerate}
\end{corollary}
\begin{proof}
\Cref{prop:collar} identifies \(X_n\) with \(E\) up to homotopy,
so the algebra embeddings and the assertions about the fundamental groups follow from
\Cref{prop:transfer-splitting}.
For the $\pi_3$ assertions, let $p\colon\operatorname{Spin}(n)\to\SO(n)$ be the spin covering. By \cite[Proposition~4.1]{Hatcher}, it induces an isomorphism on $\pi_3$. Thus \Cref{prop:transfer-splitting} implies that the composite $o_f\circ p$ induces a split injection on $\pi_3$ into $X_n$.

By the identifications in \cite{LawsonMichelsohn},
\[
 \operatorname{Spin}(3)\cong S^3,
 \qquad \operatorname{Spin}(4)\cong S^3\times S^3,
\]
so $\pi_3(\operatorname{Spin}(3))\cong\Z$ and $\pi_3(\operatorname{Spin}(4))\cong\Z^2$. The inclusions $X_n\hookrightarrow E$ are the homotopy equivalences noted above, so the same split injections hold with target $E$. Since $\pi_3(E)$ is abelian, these give the stated direct summands.

Finally, the algebra embeddings combined with the universal coefficient theorem over \(\Q\) give the claims about  Betti numbers.
\end{proof}

Now consider the Robinson sextic \cite{Robinson}
\[
 R_0=x^6+y^6+z^6-\sum_{\mathrm{sym}}x^4y^2+3x^2y^2z^2.
\]
Its ten projective zeros are the four points \([1:\pm1:\pm1]\) and the six
points obtained by choosing the zero coordinate and the relative sign of the
other two coordinates. The \(10\times10\) evaluation matrix on cubic
monomials at these representatives has determinant of absolute value
\(128\). Thus \(I_3(Z(R_0))=0\), and \(r_0=R_0/L(R_0)\in D\).
Let \(O\subseteq\SO(3)\) be the rotational octahedral group. It preserves
\(r_0\) and acts faithfully on its zeros. Since \(r_0\in D\), averaging a
positive definite matrix in \(C_{r_0}\) over \(O\) gives an invariant
\(B\in C_{r_0}\) with \(B>0\). For sufficiently small \(t>0\), set
\[
 Q=tB^{-1},\qquad f_0=\frac{r_0+\Gamma(Q)}{L(r_0+\Gamma(Q))}.
\]
As \(t\) goes to zero, \(f_0\) converges to \(r_0\notin S\), so \(f_0\notin S\) for small
\(t\); and \(Q>0\) makes \(f_0\) strictly positive. Hence \(f_0\in X_3\).
Since \(Q^{-1}=t^{-1}B\in C_{r_0}\), \eqref{eq:kkt} gives
\(R(f_0)=r_0\) after normalization. The \(O\)-fixed form \(f_0\) defines
\[
 j_3\colon\SO(3)/O\to X_3,\qquad j_3(gO)=g\cdot f_0.
\]
The composite \(\operatorname{Spin}(3)\to\SO(3)\to\SO(3)/O\) is the
universal cover, with deck group the binary octahedral group \(2O\) of order
\(48\) \cite{CoxeterBinary}. Thus \(\pi_1(\SO(3)/O)\cong2O\).
The ten-sheeted cover \(p_1\colon\mathcal C_1(X_3)\to X_3\) defines a
monodromy
\[
 \mu_3\colon\pi_1(X_3,f_0)\to\mathfrak S(Z(r_0));
\]
write \(\rho_3\colon O\hookrightarrow\mathfrak S(Z(r_0))\) for the faithful
action on the zeros.

For the quaternary case, take the Lax--Lax quartic \cite{LaxLax}
\[
 A(x)=\sum_{i=1}^5\prod_{j\ne i}(x_i-x_j).
\]
It is nonnegative and translation-invariant, so its restriction to the
sum-zero hyperplane \(W=\{y\in\R^5\mid\sum_i y_i=0\}\cong\R^4\) gives
\[
 8A(y)=L_0(y),\qquad
 L_0(y)=30\sum_{i=1}^5 y_i^4-7\left(\sum_{i=1}^5y_i^2\right)^2.
\]
Its ten projective zeros are represented by the permutations of
\((3,3,-2,-2,-2)\). After eliminating \(y_5\), the \(10\times10\) evaluation
matrix on quadratic monomials has determinant of absolute value
\(3\cdot5^{15}\). Hence \(I_2(Z(L_0))=0\) and
\(\ell_0=L_0/L(L_0)\in D\).
Let \(I=A_5\subseteq\SO(W)\) act by even coordinate permutations. The zeros
are indexed by the two-element subsets of \(\{1,2,3,4,5\}\), so this action
is faithful; also \(\ell_0\) is \(I\)-invariant. Averaging a positive
definite element of \(C_{\ell_0}\) over \(I\) gives an invariant matrix
\(B>0\). For sufficiently small \(t>0\), put
\[
 Q=tB^{-1},\qquad
 f_0=\frac{\ell_0+\Gamma(Q)}{L(\ell_0+\Gamma(Q))}.
\]
The same argument as above gives \(f_0\in X_4\) with
\(R(f_0)=\ell_0\). Thus \(f_0\) is \(I\)-fixed and defines
\(j_4\colon\SO(4)/I\to X_4\). The composite
\(\operatorname{Spin}(4)\to\SO(4)\to\SO(4)/I\) is the universal cover.
To identify its deck group, note that the double transposition \((12)(34)\)
acts on \(W\) with eigenvalues \(1,1,-1,-1\), hence as a half-turn. Its
lifts are \(\pm e e\prime\) for an orthonormal basis \(e,e\prime\) of the
rotation plane, and \((e e\prime)^2=-1\). Thus the preimage of \(I\) is a
nonsplit central double cover of \(A_5\). Since the Schur multiplier of
\(A_5\) is \(\mathbb Z/2\), its unique nontrivial central double cover is
the binary icosahedral group \(2I\) of order \(120\)
\cite{ConwaySmith,CoxeterBinary}; hence
\(\pi_1(\SO(4)/I)\cong2I\). The cover
\(p_1\colon\mathcal C_1(X_4)\to X_4\) similarly defines
\(\mu_4\colon\pi_1(X_4,f_0)\to\mathfrak S(Z(\ell_0))\); write
\(\rho_4\colon I\hookrightarrow\mathfrak S(Z(\ell_0))\) for the faithful
action.

We now prove the only missing statement from the above main theorems.

\begin{proposition}\label{prop:sextic-fundamental-group}
The maps \((j_3)_*\) and \((j_4)_*\) embed \(2O\) and \(2I\) in
\(\pi_1(X_3)\) and \(\pi_1(X_4)\), respectively.
\end{proposition}
\begin{proof}
Let \(\pi_3\colon2O\to O\), \(\pi_4\colon2I\to I\) be the covering maps,
and let \(q_3\colon\SO(3)\to\SO(3)/O\) and
\(q_4\colon\SO(4)\to\SO(4)/I\) be the quotient maps. Equivariance of \(R\)
and path lifting give
\[
 \mu_3\circ(j_3)_*=\rho_3\circ\pi_3,\qquad
 \mu_4\circ(j_4)_*=\rho_4\circ\pi_4.
\]
Faithfulness implies \(\ker(j_3)_*\subseteq\ker\pi_3=\{\pm1\}\) and
\(\ker(j_4)_*\subseteq\ker\pi_4=\{\pm1\}\). Moreover,
\[
 (q_3)_*\pi_1(\SO(3))=\ker\pi_3,\qquad
 (q_4)_*\pi_1(\SO(4))=\ker\pi_4,
\]
while \((j_3\circ q_3)_*=(o_{f_0})_*\) and
\((j_4\circ q_4)_*=(o_{f_0})_*\) are injective by
\Cref{prop:orbit-pi1}. Thus \((j_3)_*\) and \((j_4)_*\) are injective
on \(\ker\pi_3\) and \(\ker\pi_4\), respectively. The kernel inclusions
above then show that both maps are injective.
\end{proof}

\section*{Declaration on the use of generative AI}
The authors used  artificial intelligence systems for exploring ideas, streamlining proofs and improving the presentation of the paper. All AI-assisted suggestions were reviewed and edited by the authors, who take responsibility for the final manuscript.

\end{document}